\documentclass{amsart}

\usepackage{amsmath}
\usepackage{amsfonts}
\usepackage{amsthm} 
\usepackage{amssymb}

\usepackage{graphicx} 

\usepackage{color} 

\newtheorem{thm}{Theorem}[section]

\newtheorem{cor}[thm]{Corollary}

\numberwithin{equation}{section} 

\makeatletter 
\@namedef{subjclassname@2010}{2010 Mathematics Subject Classification} 
\makeatother 

\begin{document}

\title[On oscillatory integrals with degenerate phase functions]
{On asymptotic expansions of oscillatory integrals with phase functions expressed by a product of positive real power function and real analytic function in one variable} 

\author{Toshio NAGANO} 

\address{Department of Liberal Arts, Faculty of Science and Technology, Tokyo University of Science, 2641, Yamazaki, Noda, Chiba 278-8510, JAPAN} 
\email{tonagan@rs.tus.ac.jp}

\thanks{
The author was supported by Tokyo University of Science Graduate School doctoral program scholarship and an exemption of the cost of equipment from 2016 to 2018 and would like to thank to Professor Minoru Ito for giving me an opportunity of studies and preparing the environment.
} 



\date {October 22, 2020} 

\begin{abstract} 
In this paper, by using asymptotic expansions of oscillatory integrals with positive real power phase functions in one variable, 
we obtain asymptotic expansions of oscillatory integrals with phase functions expressed by a product of positive real power function and real analytic function in one variable. 
Moreover we show an example which we can compute all coefficients of terms in asymptotic expansions concretely. 
\end{abstract} 

\maketitle 

\section{Introduction} 

We study 
for {\bf degenerate} phase functions $\phi(x)$ and suitable amplitude functions $a(x)$ in one variable, 
asymptotic expansions of oscillatory integrals: 
\begin{align} 
Os\text{-}\int_{-\infty}^{\infty} e^{i \lambda \phi(x)} a (x) dx 
:= \lim_{\varepsilon \to +0} \int_{-\infty}^{\infty} 
e^{i \lambda \phi(x)} a (x) \chi (\varepsilon x) dx, 
\notag 
\end{align} 
as $\lambda \to \infty$, 
where $\chi$ is a rapidly decreasing function of class $C^{\infty}$ on $\mathbb{R}$ with $\chi(0) = 1$ and $0 < \varepsilon < 1$. 

In the previous paper \cite{Nagano-Miyazaki02}, 
we studied the case of $\phi(x) = x^{p}$ for any positive real number $p>0$ and $a(x) \in \mathcal{A}^{\tau}_{\delta}(\mathbb{R})$, where for any $\tau \in \mathbb{R}$ and $\delta \in [-1,p-1)$, $\mathcal{A}^{\tau}_{\delta}(\mathbb{R})$ is the set of all $a \in C^{\infty}(\mathbb{R})$ such that for any $l \in \mathbb{Z}_{\geq 0}$, 
\begin{align} 
|a|^{(\tau)}_{l} 
:= \max_{k \leq l} \sup_{x \in \mathbb{R}} (1+|x|^{2})^{-(\tau + \delta k)/2} |a^{(k)}(x)| < \infty. 
\notag 
\end{align} 
Then we obtained the following results (Theorem 3.3 (i) and Theorem 5.2 in \cite{Nagano-Miyazaki02}): 

\begin{thm} 
\label{Lax02} 
Assume that $\lambda > 0$ and $p > 0$. 
Let $\varphi \in C^{\infty}_{0}(\mathbb{R})$. 
Then there exist the following oscillatory integrals, and the following holds: 
\begin{align} 
Os\text{-}\int_{0}^{\infty} e^{\pm i\lambda x^{p}} \varphi(x) dx 
= \int_{0}^{\infty} e^{\pm i\lambda x^{p}} \varphi(x) dx, 
\notag 
\end{align} 
where double signs $\pm$ are in the same order. 
\end{thm} 

\begin{thm} 
\label{th02} 
Assume that $\lambda > 0$ and $a \in \mathcal{A}^{\tau}_{\delta}(\mathbb{R})$. 
Then the following hold: 
 \begin{enumerate} 
\item[(i)] 
If $p>0$, 
then for any $N \in \mathbb{N}$ such that $N \geq p+1$, as $\lambda \to \infty$, 
\begin{align} 
Os\text{-}\int_{0}^{\infty} e^{\pm i\lambda x^{p}} a(x) dx 
&= \sum_{k=0}^{N-[p]-1} \tilde{I}_{p,k+1}^{\pm} \frac{a^{(k)}(0)}{k!} \lambda ^{-\frac{k+1}{p}} + O\left( \lambda^{-\frac{N-p+1}{p}} \right), 
\notag 
\end{align} 
where $[p]$ is the Gauss' symbol, that is, $[p] \in \mathbb{Z}$ such that $p-1 < [p] \leq p$, 
\begin{align} 
\tilde{I}_{p,k+1}^{\pm} 
= p^{-1} e^{\pm i\frac{\pi}{2} \frac{k+1}{p}} \varGamma \left( \frac{k+1}{p} \right), 
\label{tilde_I_p_k+1_pm} 
\end{align} 
and double signs $\pm$ are in the same order. 
 
\item[(ii)] 
If $p=m \in \mathbb{N}$, 
then for any $N \in \mathbb{N}$ such that $N > m$, as $\lambda \to \infty$, 
\begin{align} 
Os\text{-}\int_{0}^{\infty} e^{\pm i\lambda (-y)^{m}} a(-y) dy 
&= \sum_{k=0}^{N-m-1} (-1)^{k} \tilde{I}_{m,k+1}^{\pm \pm^{m}} \frac{a^{(k)}(0)}{k!} \lambda ^{-\frac{k+1}{m}} + O\left( \lambda^{-\frac{N-m+1}{m}} \right), 
\notag 
\end{align} 
where 
\begin{align} 
\tilde{I}_{m,k+1}^{\pm \pm^{m}} 
&= m^{-1} e^{\pm (-1)^{m} i\frac{\pi}{2} \frac{k+1}{m}} \varGamma \left( \frac{k+1}{m} \right), 
\label{tilde_I_m_k+1_pm_pm_m} 
\end{align} 
and double signs $\pm$ are in the same order. 
\end{enumerate} 
\end{thm} 

In this paper, by using Theorem \ref{Lax02}, as an application of Theorem \ref{th02}, 
we obtain the following results: 
\begin{thm} 
\label{main theorem} 
Assume that $\lambda > 0$ and $p > 0$. 
Let $\{ a_{j} \}_{j \in \mathbb{N}}$ be a sequence of real numbers 
such that $l_{0} := \sup_{j \in \mathbb{N}} |a_{j}|^{1/j} < \infty$, 
and let 
\begin{align} 
\phi (x) 
:= x^{p} \bigg( 1 + \sum_{j=1}^{\infty} a_{j} x^{j} \bigg) 
\notag 
\end{align} 
for $|x| < R_{0} := (2l_{0})^{-1}$ where if $l_{0}=0$, then $R_{0}=\infty$. 
Then 
there exists a neighborhood $U$ of the origin in $\mathbb{R}$ such that $U \subset (-R_{0},R_{0})$, 
and a diffeomorphism $x = \varPhi (y)$ for $x,y \in U$ such that $x (1 + \sum_{j=1}^{\infty} a_{j} x^{j}) ^{1/p} = y$, $\varPhi (U \cap [0,\infty)) = U \cap [0,\infty)$ and $\varPhi (U \cap (-\infty,0]) = U \cap (-\infty,0]$, 
for any $a \in C^{\infty}_{0}(\mathbb{R})$ such that $\mathrm{supp}~a \subset U$, 
the following hold: 

\begin{enumerate} 
\item[(i)] 
For any $N \in \mathbb{N}$ such that $N \geq p+1$, as $\lambda \to \infty$, 
\begin{align} 
&\int_{0}^{\infty} e^{\pm i\lambda x^{p} \left( 1 + \sum_{j=1}^{\infty} a_{j} x^{j} \right)} a(x) dx \notag \\ 
&= \sum_{k=0}^{N-[p]-1} 
\frac{\tilde{I}_{p,k+1}^{\pm}}{k!} \left( \frac{d}{dy} \right)^{k} \bigg|_{y=0} \left\{ a(\varPhi (y)) \frac{d\varPhi}{dy}(y) \right\} \lambda ^{-\frac{k+1}{p}} 
+ O\left( \lambda ^{-\frac{N-p+1}{p}} \right), 
\notag 
\end{align} 
where $\tilde{I}_{p,k+1}^{\pm}$ is given by \eqref{tilde_I_p_k+1_pm} and double signs $\pm$ are in the same order. 

\item[(ii)] 
If $p = m \in \mathbb{N}$, 
then for any $N \in \mathbb{N}$ such that $N > m$, as $\lambda \to \infty$, 
\begin{align} 
&\int_{-\infty}^{\infty} e^{\pm i\lambda x^{m} \left( 1 + \sum_{j=1}^{\infty} a_{j} x^{j} \right)} a(x) dx \notag \\ 
&= \sum_{k=0}^{N-m-1} 
\frac{c_{m,k}^{\pm}}{k!} \left( \frac{d}{dy} \right)^{k} \bigg|_{y=0} \left\{ a(\varPhi (y)) \frac{d\varPhi}{dy}(y) \right\} \lambda ^{-\frac{k+1}{m}} + O\left( \lambda ^{-\frac{N-m+1}{m}} \right), 
\notag 
\end{align} 
where 
\begin{align} 
c_{m,k}^{\pm} 
= m^{-1} \left\{ e^{\pm i\frac{\pi}{2} \frac{k+1}{m}} + (-1)^{k+1} e^{\pm (-1)^{m} i\frac{\pi}{2} \frac{k+1}{m}} \right\} \varGamma \left( \frac{k+1}{m} \right), 
\notag 
\end{align} 
and double signs $\pm$ are in the same order. 
\end{enumerate} 
\end{thm} 

In particular, when $a_{j}=1/j!$, by using differential coefficients of all orders, of Lambert W-function $X=W(Y)$ (cf.\cite{CGHJK}) at $Y=0$, which is an inverse function of $Y=Xe^{X}$ defined on a certain neighborhood of the origin in $\mathbb{R}$, then we can compute all coefficients of terms in asymptotic expansions concretely as follows: 
\begin{cor} 
\begin{enumerate} 
\item[(i)] 
\begin{align} 
&\int_{0}^{\infty} e^{\pm i\lambda x^{p} e^{x}} a(x) dx \notag \\ 
&= \sum_{k=0}^{N-[p]-1} 
\frac{\tilde{I}_{p,k+1}^{\pm}}{k!} \left( \frac{d}{dy} \right)^{k} \bigg|_{y=0} \left\{ a(\varPhi (y)) \frac{dW}{dY} \bigg( \frac{y}{p} \bigg) \right\} \lambda ^{-\frac{k+1}{p}} 
+ O\left( \lambda ^{-\frac{N-p+1}{p}} \right), 
\notag 
\end{align} 
where $\tilde{I}_{p,k+1}^{\pm}$ is given by \eqref{tilde_I_p_k+1_pm} and double signs $\pm$ are in the same order. 

\item[(ii)] 
If $p = m \in \mathbb{N}$, 
then for any $N \in \mathbb{N}$ such that $N > m$, as $\lambda \to \infty$, 
\begin{align} 
&\int_{-\infty}^{\infty} e^{\pm i\lambda x^{m} e^{x}} a(x) dx \notag \\ 
&= \sum_{k=0}^{N-m-1} 
\frac{c_{m,k}^{\pm}}{k!} \left( \frac{d}{dy} \right)^{k} \bigg|_{y=0} \left\{ a(\varPhi (y)) \frac{dW}{dY} \bigg( \frac{y}{m} \bigg) \right\} \lambda ^{-\frac{k+1}{m}} + O\left( \lambda ^{-\frac{N-m+1}{m}} \right), 
\notag 
\end{align} 
where 
\begin{align} 
c_{m,k}^{\pm} 
= m^{-1} \left\{ e^{\pm i\frac{\pi}{2} \frac{k+1}{m}} + (-1)^{k+1} e^{\pm (-1)^{m} i\frac{\pi}{2} \frac{k+1}{m}} \right\} \varGamma \left( \frac{k+1}{m} \right), 
\notag 
\end{align} 
and double signs $\pm$ are in the same order. 
\end{enumerate} 
\end{cor} 

To the end of \S1, we note notation used in this paper: 

$C^{\infty}(\mathbb{R}^{n})$ is 
the set of complex-valued functions of class $C^{\infty}$ on $\mathbb{R}^{n}$. 
$C^{\infty}_{0}(\mathbb{R}^{n})$ is 
the set of all $f \in C^{\infty}(\mathbb{R}^{n})$ with compact support. 

$O$ means the Landau's symbol, that is, 
$f(x) = O(g(x))~(x \to a)$ if $|f(x)/g(x)|$ 
is bounded as $x \to a$ for functions $f$ and $g$, where $a \in \mathbb{R} \cup \{ \pm \infty \}$. 

\section{Proofs of main results} 

\begin{proof}[Proof of Theorem 1.3]  
(i) 
Since $l := \mathrm{lim~sup}_{j\to \infty} |a_{j}|^{1/j} \leq \sup_{j \in \mathbb{N}} |a_{j}|^{1/j} =: l_{0} < \infty$, 
by Cauchy-Hadamard's formula, 
the radius of convergence of the power series $\sum_{j=1}^{\infty} a_{j} x^{j}$ is $R := 1/l >0$. 
Moreover 
if $|x| < R_{0} := (2l_{0})^{-1} \leq (2l)^{-1} < R$, 
then 
\begin{align} 
\bigg| \sum_{j=1}^{\infty} a_{j} x^{j} \bigg| 
\leq \sum_{j=1}^{\infty} |a_{j}| |x|^{j} 
< \sum_{j=1}^{\infty} |a_{j}| \Big\{ \big( 2 \sup_{j \in \mathbb{N}} |a_{j}|^{1/j} \big)^{-1} \Big\}^{j} 
< \sum_{j=1}^{\infty} 2^{-j} 
= 1. 
\notag 
\end{align} 
Hence 
if $|x| < R_{0}$, 
since $1 + \sum_{j=1}^{\infty} a_{j} x^{j} > 0$, 
then we can define 
\begin{align} 
f(x) 
:= x \bigg( 1 + \sum_{j=1}^{\infty} a_{j} x^{j} \bigg) ^{1/p} 
\notag 
\end{align} 
for $x \in (-R_{0},R_{0})$. 
Since $f$ is a function of class $C^{\infty}$ on $(-R_{0},R_{0})$ with $f(0) = 0$ and $f'(0) = 1$, 
by inverse function theorem, 
there exists a neighborhood $U$ of the origin in $\mathbb{R}$ such that $U \subset (-R_{0},R_{0})$, and a diffeomorphism $x = \varPhi (y)$ for $x,y \in U$ such that $f(x) = y$, $\varPhi (U \cap [0,\infty)) = U \cap [0,\infty)$ and $\varPhi (U \cap (-\infty,0]) = U \cap (-\infty,0]$. 
Then for any $a \in C^{\infty}_{0}(\mathbb{R})$ such that $\mathrm{supp}~a \subset U$, 
the following improper integrals are absolutely convergent: 
\begin{align} 
\int_{0}^{\infty} e^{\pm i\lambda x^{p} \left( 1 + \sum_{j=1}^{\infty} a_{j} x^{j} \right)} a(x) dx~\text{and}~
\int_{-\infty}^{0} e^{\pm i\lambda x^{m} \left( 1 + \sum_{j=1}^{\infty} a_{j} x^{j} \right)} a(x) dx~\text{for $m \in \mathbb{N}$}. 
\notag 
\end{align} 
By change of variable $x = \varPhi (y)$ on $U \cap [0,\infty)$, 
since $x^{p} (1 + \sum_{j=1}^{\infty} a_{j} x^{j}) = (f(x))^{p} = y^{p}$, 
\begin{align} 
\int_{0}^{\infty} e^{\pm i\lambda x^{p} \left( 1 + \sum_{j=1}^{\infty} a_{j} x^{j} \right)} a(x) dx 
= \int_{0}^{\infty} e^{\pm i\lambda y^{p}} a(\varPhi (y)) \frac{d\varPhi}{dy}(y) dy. 
\notag 
\end{align} 
Here since $a(\varPhi (y)) \frac{d\varPhi}{dy}(y) \in C^{\infty}_{0}(\mathbb{R})$, 
by Theorem \ref{Lax02} and Theorem \ref{th02} (i), 
for any $N \in \mathbb{N}$ such that $N \geq p+1$, as $\lambda \to \infty$, 
\begin{align} 
&\int_{0}^{\infty} e^{\pm i\lambda x^{p} \left( 1 + \sum_{j=1}^{\infty} a_{j} x^{j} \right)} a(x) dx 
= Os\text{-}\int_{0}^{\infty} e^{\pm i\lambda y^{p}} a(\varPhi (y)) \frac{d\varPhi}{dy}(y) dy \notag \\ 
&= \sum_{k=0}^{N-[p]-1} 
\frac{\tilde{I}_{p,k+1}^{\pm}}{k!} \left( \frac{d}{dy} \right)^{k} \bigg|_{y=0} \left\{ a(\varPhi (y)) \frac{d\varPhi}{dy}(y) \right\} \lambda ^{-\frac{k+1}{p}} 
+ O\left( \lambda ^{-\frac{N-p+1}{p}} \right). 
\notag 
\end{align} 

(ii) 
If $p = m \in \mathbb{N}$, then we take $N \in \mathbb{N}$ such that $N > m$, since $N \geq m+1$ and $[m]=m$, by (i), as $\lambda \to \infty$, 
\begin{align} 
&\int_{0}^{\infty} e^{\pm i\lambda x^{m} \left( 1 + \sum_{j=1}^{\infty} a_{j} x^{j} \right)} a(x) dx \notag \\ 
&= \sum_{k=0}^{N-m-1} 
\frac{\tilde{I}_{m,k+1}^{\pm}}{k!} \left( \frac{d}{dy} \right)^{k} \bigg|_{y=0} \left\{ a(\varPhi (y)) \frac{d\varPhi}{dy}(y) \right\} \lambda ^{-\frac{k+1}{m}} 
+ O\left( \lambda ^{-\frac{N-m+1}{m}} \right), 
\label{int_0_infty} 
\end{align} 
where $\tilde{I}_{m,k+1}^{\pm}$ is given by \eqref{tilde_I_p_k+1_pm} and double signs $\pm$ are in the same order. 

On the other hand, by (i), the following improper integrals is absolutely convergent, 
and by change of variable $x = \varPhi (y)$ on $U \cap (-\infty,0]$ and $y=-u$, 
since $x^{m} (1 + \sum_{j=1}^{\infty} a_{j} x^{j}) = (f(x))^{m} = y^{m} = (-u)^{m}$, the following holds: 
\begin{align} 
\int_{-\infty}^{0} e^{\pm i\lambda x^{m} \left( 1 + \sum_{j=1}^{\infty} a_{j} x^{j} \right)} a(x) dx 
&= \int_{-\infty}^{0} e^{\pm i\lambda y^{m}} a(\varPhi (y)) \frac{d\varPhi}{dy}(y) dy \notag \\ 
&= -\int_{0}^{\infty} e^{\pm i\lambda (-u)^{m}} a(\varPhi(-u)) \frac{d\varPhi}{du}(-u) du. 
\notag 
\end{align} 
Here 
since $a(\varPhi (-u)) \frac{d\varPhi}{du}(-u) \in C^{\infty}_{0}(\mathbb{R})$, 
by Theorem \ref{Lax02} and Theorem \ref{th02} (ii), 
for any $N \in \mathbb{N}$ such that $N > m$, as $\lambda \to \infty$, 
\begin{align} 
&\int_{-\infty}^{0} e^{\pm i\lambda x^{m} \left( 1 + \sum_{j=1}^{\infty} a_{j} x^{j} \right)} a(x) dx 
= -Os\text{-}\int_{0}^{\infty} e^{\pm i\lambda (-u)^{m}} a(\varPhi(-u)) \frac{d\varPhi}{du}(-u) du \notag \\ 
&= -\sum_{k=0}^{N-m-1} 
\frac{\tilde{I}_{m,k+1}^{\pm \pm^{m}}}{k!} \left( \frac{d}{du} \right)^{k} \bigg|_{u=0} \left\{ a(\varPhi(-u)) \frac{d\varPhi}{du}(-u) \right\} \lambda ^{-\frac{k+1}{m}} 
+ O\left( \lambda ^{-\frac{N-m+1}{m}} \right), 
\label{int_-infty_0} 
\end{align} 
where $\tilde{I}_{m,k+1}^{\pm \pm^{m}}$ is given by \eqref{tilde_I_m_k+1_pm_pm_m} and double signs $\pm$ are in the same order. 

Therefore by \eqref{int_0_infty} and \eqref{int_-infty_0}, 
for any $N \in \mathbb{N}$ such that $N > m$, as $\lambda \to \infty$, 
\begin{align} 
&\int_{-\infty}^{\infty} e^{\pm i\lambda x^{m} \left( 1 + \sum_{j=1}^{\infty} a_{j} x^{j} \right)} a(x) dx \notag \\ 
&= \sum_{k=0}^{N-m-1} 
\frac{c_{m,k}^{\pm}}{k!} \left( \frac{d}{dy} \right)^{k} \bigg|_{y=0} \left\{ a(\varPhi (y)) \frac{d\varPhi}{dy}(y) \right\} \lambda ^{-\frac{k+1}{m}} + O\left( \lambda ^{-\frac{N-m+1}{m}} \right), 
\notag 
\end{align} 
where 
\begin{align} 
c_{m,k}^{\pm} 
&=\tilde{I}_{m,k+1}^{\pm} + (-1)^{k+1} \tilde{I}_{m,k+1}^{\pm \pm^{m}} \notag \\ 
&= m^{-1} \left\{ e^{\pm i\frac{\pi}{2} \frac{k+1}{m}} + (-1)^{k+1} e^{\pm (-1)^{m} i\frac{\pi}{2} \frac{k+1}{m}} \right\} \varGamma \left( \frac{k+1}{m} \right), 
\notag 
\end{align} 
and double signs $\pm$ are in the same order. 
\end{proof} 

\begin{proof}[Proof of Corollary 1.4]  
(i) In Theorem 1.3 (i), if $a_{j}=1/j!$, since 
\begin{align} 
\sup_{j \in \mathbb{N}} |1/j!|^{1/j} = \sup_{j \in \mathbb{N}} \prod_{k=1}^{j} (1/k)^{1/j} \leq 1, 
\notag 
\end{align} 
then 
\begin{align} 
f(x) 
= x \bigg( 1 + \sum_{j=1}^{\infty} \frac{x^{j}}{j!} \bigg) ^{1/p} 
= x e^{x/p}. 
\notag 
\end{align} 
Here let $y=f(x)$, since $y/p = (x/p) e^{x/p}$, then
\begin{align} 
\frac{x}{p} = W \bigg( \frac{y}{p} \bigg), 
\notag 
\end{align} 
where $X=W(Y)$ is a Lambert W-function on $(1/p) U \cap [0,\infty)$. 
And then since $\varPhi(y) = x = pW(y/p)$, 
\begin{align} 
\frac{d\varPhi}{dy}(y) 
= \frac{d}{dy} \bigg\{ p W \bigg( \frac{y}{p} \bigg) \bigg\} 
= \frac{dW}{dY} \bigg( \frac{y}{p} \bigg). 
\notag 
\end{align} 

 (ii) When $p=m \in \mathbb{N}$, if $x < 0$, then we take a Lambert W-function $X=W(Y)$ on $(1/m) U \cap (\infty,0]$, 
we can prove (ii) in the same manner to (i). 
\end{proof}

\end{document}